\documentclass[journal]{IEEEtran}

 \usepackage{graphicx} 
  \DeclareGraphicsExtensions{.eps}
  
  \usepackage{tikz}
  \usepackage{pgfplots,pgfplotstable}
  \usetikzlibrary{pgfplots.groupplots}
  \pgfplotsset{compat=1.5}

\usepackage{multirow}

\usepackage{algorithm,algorithmic}

\usepackage{enumitem} 

\usepackage{subcaption} 

\usepackage[font=footnotesize]{caption} 

\usepackage{amsmath}

\usepackage{textcomp}
\usetikzlibrary{shapes,arrows}

\usepackage{amsthm} 
\newtheorem{prop}{Proposition}

\usepackage{empheq} 

\usepackage{mathrsfs} 
\newsavebox\foobox
\newlength{\foodim}
\newcommand{\slantbox}[2][0]{\mbox{%
        \sbox{\foobox}{#2}%
        \foodim=#1\wd\foobox
        \hskip \wd\foobox
        \hskip -0.5\foodim
        \pdfsave
        \pdfsetmatrix{1 0 #1 1}%
        \llap{\usebox{\foobox}}%
        \pdfrestore
        \hskip 0.5\foodim
}}
\def\Laplace{\slantbox[-.45]{$\mathscr{L}$}}

\usepackage{stfloats} 

\usepackage[hyphens]{url}
\usepackage[hidelinks]{hyperref}
\hypersetup{colorlinks=false,linkcolor=black,citecolor=black,filecolor=black,urlcolor=black,breaklinks=true}
\usepackage[normalem]{ulem} 

\usepackage{cite} 

\usepackage{color}

\usepackage{cancel} 

\usepackage{array} 

\usepackage[export]{adjustbox} 

\begin{document}

\bstctlcite{IEEEexample:BSTcontrol}

\title{Conditions for Regional Frequency Stability in Power System Scheduling---Part I: Theory}

%
\author{Luis~Badesa,~\IEEEmembership{Member,~IEEE,}
        Fei~Teng,~\IEEEmembership{Member,~IEEE,}
        and~Goran~Strbac,~\IEEEmembership{Member,~IEEE}
\thanks{A portion of this work has been supported by the UK Engineering and Physical Sciences Research Council under project `Integrated Development of Low-Carbon Energy Systems’ (IDLES), grant EP/R045518/1.}
\thanks{
The authors are with the Department of Electrical and Electronic Engineering, Imperial College London, SW7 2AZ London, U.K. (email: luis.badesa@imperial.ac.uk, f.teng@imperial.ac.uk, g.strbac@imperial.ac.uk).}
%
}

%
%

\markboth{IEEE Transactions on Power Systems, April~2021}%
{Shell \MakeLowercase{\textit{et al.}}: Bare Demo of IEEEtran.cls for IEEE Journals}
%



\maketitle

\begin{abstract}
This paper considers the phenomenon of distinct regional frequencies recently observed in some power systems. First, a reduced-order mathematical model describing this behaviour is developed. Then, techniques to solve the model are discussed, demonstrating that the post-fault frequency evolution in any given region is equal to the frequency evolution of the Centre Of Inertia plus certain inter-area oscillations. This finding leads to the deduction of conditions for guaranteeing frequency stability in all regions of a power system, a deduction performed using a mixed analytical-numerical approach that combines mathematical analysis with regression methods on simulation samples. The proposed stability conditions are linear inequalities that can be implemented in any optimisation routine allowing the co-optimisation of all existing ancillary services for frequency support: inertia, multi-speed frequency response, load damping and an optimised largest power infeed. This is the first reported mathematical framework with explicit conditions to maintain frequency stability in a power system exhibiting inter-area oscillations in frequency.
\end{abstract}


%
\IEEEpeerreviewmaketitle

\section*{Nomenclature}
\addcontentsline{toc}{section}{Nomenclature}

\subsection*{Indices and Sets}
\begin{IEEEdescription}[\IEEEusemathlabelsep\IEEEsetlabelwidth{$i,\,\, j,\,\, n$}]
\setlength\itemsep{0.3em} 
\item[$i,\,\, j,\,\, n$] All-purpose indices.
\item[$\mathcal{J}_i$] Set of neighbouring regions to region $i$.
\end{IEEEdescription}

\subsection*{Constants and Parameters}
\begin{IEEEdescription}[\IEEEusemathlabelsep\IEEEsetlabelwidth{$\textrm{RoCoF}_{\textrm{max}}$}]
\setlength\itemsep{0.3em} 
\item[$\Delta f_{\textrm{max}}$] Maximum admissible frequency deviation at the nadir (Hz).
\item[$\Delta f^{\textrm{ss}}_{\textrm{max}}$] Maximum admissible frequency deviation at quasi-steady-state (Hz).
\item[$\delta_i^\textrm{ss}$] Pre-fault, steady-state phase angle of voltage in bus $i$ (rad).
\item[$\phi_i$] Phase shift of inter-area oscillations for region $i$ (rad).
\item[$\omega_i$] Angular frequency of inter-area oscillations for region $i$ (rad/s).
\item[$\textrm{a}_i$] Attenuation factor of inter-area oscillations for region $i$ ($\textrm{s}^{-1}$).
\item[$\textrm{A}_i$] Amplitude of inter-area oscillations for region $i$ (Hz).
\item[$\textrm{D}_i$] Load damping factor in region $i$ (\%/Hz).
\item[$H_i$] System inertia in region $i$ (MW$\cdot \textrm{s}$).
\item[$\textrm{P}^{\textrm{D}}_i$] Total demand in region $i$ (MW).
\item[$P^{\textrm{L}}_i$] Largest power infeed in region $i$ (MW).
\item[$R_i$] Total PFR in region $i$ (MW).
\item[$\textrm{RoCoF}_{\textrm{max}}$] Maximum admissible RoCoF (Hz/s).
\item[$\textrm{T}_{\textrm{g}}$] Delivery time of PFR (s).
\item[$\textrm{T}_{i,j}$] Electrical stiffness of the transmission line connecting buses $i$ and $j$ (MW).
\item[$\textrm{V}_i$] Voltage magnitude in bus $i$ (kV).
\item[$\textrm{X}_{i,j}$] Reactance of transmission between buses $i$ and $j$ ($\Omega$).
\end{IEEEdescription}

\subsection*{Functions and Operators}
\begin{IEEEdescription}[\IEEEusemathlabelsep\IEEEsetlabelwidth{$\Delta\textrm{P}^{\textrm{import}}_i(t)$}]
\setlength\itemsep{0.3em} 
\item[$\,\, \overline{\vphantom{Z_i} \; \, \cdot \; \,}$] Conjugate of a complex number.
\item[$\Delta f_i(t)$] Time-evolution of post-fault frequency deviation from nominal state in region $i$ (Hz).
\item[$\Delta F_i(s)$] Laplace transform of $\Delta f_i(t)$ (Hz).
\item[$\Delta\textrm{P}^{\textrm{import}}_i(t)$] Deviation from steady-state power import to region $i$, after an outage (MW).
\item[$\delta_i(t)$] Post-fault phase angle of voltage in bus $i$ (rad).
\item[$\Laplace{\{\,\cdot\,\}}$] Laplace transform operator.
\item[$\Laplace^{-1}{\{\,\cdot\,\}}$] Inverse Laplace transform operator.
\item[$\textrm{PFR}_i(t)$] Time-evolution of PFR in region $i$ (MW).
\item[$\textrm{sup}\{\,\cdot\,\}$] Supremum of a set.
\item[$t^*$] Time when the frequency nadir occurs (s).
\end{IEEEdescription}

\section{Introduction}
%
%
%
%
\IEEEPARstart{M}{aintaining} system frequency within acceptable limits is critical for the secure operation of a power grid. In the event of a generation outage, the subsequent frequency drop is contained by certain ancillary services: system inertia, load damping and Primary Frequency Response (PFR). While these ancillary services were widely available in grids dominated by thermal generators as a by-product of energy, the increasing penetration of non-synchronous Renewable Energy Sources (RES) has greatly reduced the level of system inertia, therefore increasing the risk of violating frequency stability.

Frequency stability is studied through comprehensive dynamic simulations of the power system, while closed-form stability conditions (which can be used by system operators to procure ancillary services) have been deduced in the literature from the swing equation \cite{KundurBook}. The swing equation is a simplification of the actual frequency dynamics in a power grid, which considers a single-bus representation assuming that all of the system's generators move coherently as a single lumped mass (Centre Of Inertia concept) and load damping is modelled as a single constant. While the uniform frequency model has provided a precise representation of frequency dynamics in systems dominated by thermal units, recent studies have shown that the Centre Of Inertia (COI) representation can be inaccurate in modern grids, in which RES are typically located in remote areas far from load centres, creating a non-uniform distribution of inertia. These spatial gradients for inertia cause distinct regional frequencies.

Geographical discrepancies in frequency have been observed after frequency events in recent years by utilities all over the world. Several recent publications and reports \cite{InmaMultiArea,DouglasWilsonSmartFreq,OverbyeLocationalInertia,GoranAndersson1} highlight this issue. The study in \cite{NatureDatabaseFreq} analyses measurements of electric power grid frequencies across several synchronous areas on three continents, showing how fluctuations change from local oscillations to a homogeneous bulk behaviour. Distributions of inertia within the network have been shown to affect its frequency behaviour. Reference \cite{TimeDelaysMultiFreq} demonstrates that non-uniform distributions of inertia in continental Europe can compromise frequency stability, given current control mechanism and delays. The work in \cite{AtiaInertiaHeterogeneity} performs a probabilistic frequency-stability analysis from dynamic simulations of systems with heterogeneous distributions of inertia, while \cite{DorflerOptimalPlacement} proposes a method for optimal placement of virtual inertia based on minimising deviation from nominal frequency after an outage.

Analytical approaches have been proposed to account for these frequency variations within an electricity grid. The formulation in \cite{CostasVournasMultiArea} uses a quasi-steady-state approximation of inter-area modes that provide good accuracy compared to detailed dynamic models.
The authors in \cite{FrequencyDivider} developed a simple algebraic formula for accurately estimating frequency in any of the grid's buses, which can significantly decrease the computation of dynamic simulations. However, to date no work has deduced conditions for limiting the frequency drop within pre-defined thresholds, key for constraining a scheduling problem that considers ancillary services. 

Considering the uniform-frequency model in systems with non-uniform inertia distribution has some dangers: the actual need for frequency ancillary services would be underestimated, leading to higher regional Rate of Change of Frequency (RoCoF) and lower regional nadirs than expected. If the generation scheduling or ancillary-services market is not appropriately constrained to reflect the distinct regional frequencies, unexpected tripping of RoCoF relays and triggering of Under-Frequency Load Shedding (UFLS) could take place, which in turn could cause cascading outages potentially leading to a blackout.

In this context, the present paper focuses on generalising the currently available conditions for guaranteeing frequency stability to account for the spatial variations in post-fault frequency dynamics. In order to do so, mathematical inequalities representing the stability boundary are deduced from a spatial swing model, which generalises the swing equation by considering $N$ different regions in the grid coupled by ac transmission lines. By solving this spatial swing model, constraints for guaranteeing regional frequency stability can be obtained, to be later implemented in optimisation routines. Nevertheless, solving this model analytically is a challenging task, even impossible for the simplest systems as will be demonstrated in coming sections of this paper. To overcome this difficulty, in this work we propose a mixed analytical-numerical method to obtain the constraints, using theoretical mathematical techniques combined with regression methods on simulation samples. 

The key contributions of this work are:
\begin{enumerate}
	\item To the best of our knowledge, this is the first work to deduce closed-form conditions for post-fault regional frequency stability in a power grid, securing the system against oscillations in frequency that exacerbate the Centre of Inertia frequency drop.
	
	\item The proposed constraints allow the co-optimisation of all existing frequency services, including inertia, multi-speed frequency response, load damping and an optimised largest power infeed. An intuitive formulation for the constraints is provided, which are suitable to be implemented in optimisation problems like market clearing and Unit Commitment.
	
\end{enumerate}

This paper is organised as follows: Section~\ref{SectionDynamicSimulations} analyses the causes of distinct regional frequencies, while the dynamic model and a discussion on how to solve it is presented in Section~\ref{sectionDynamicModelMultiRegion}. Sections~\ref{Conditions2Region} and \ref{ConditionsNregions} introduce the proposed stability conditions for two-region and $N$-region systems. Finally, Section~\ref{sectionConclusion} gives the conclusion.

\section{Causes of distinct regional frequencies} \label{SectionDynamicSimulations}

Before undertaking the deduction of mathematical constraints that would guarantee regional frequency security, we analyse here some causes of distinct post-fault frequencies in a power system. Let's consider the Great Britain (GB) system split in two areas, roughly corresponding to Scotland and England: Scotland has significant wind resources while most load centres are located in England, a fact that drives a non-uniform inertia distribution.

The two-region England-Scotland system was simulated in MATLAB/Simulink, using a mathematical model that will be discussed in Section \ref{sectionDynamicModelMultiRegion}. Simulations were run by changing the system operating condition, i.e.~changing the inertia, load damping and frequency response in each region, the location of the generation outage (placed alternatively in Scotland and England) and the strength of the electrical interconnection between the regions (driven by the impedance of the connecting transmission corridors).

\begin{figure}[!t]
	\centering
    \begin{subfigure}[b]{0.5\textwidth}
        \hspace*{-8mm}
        \centering
        \includegraphics[width=3.1in]{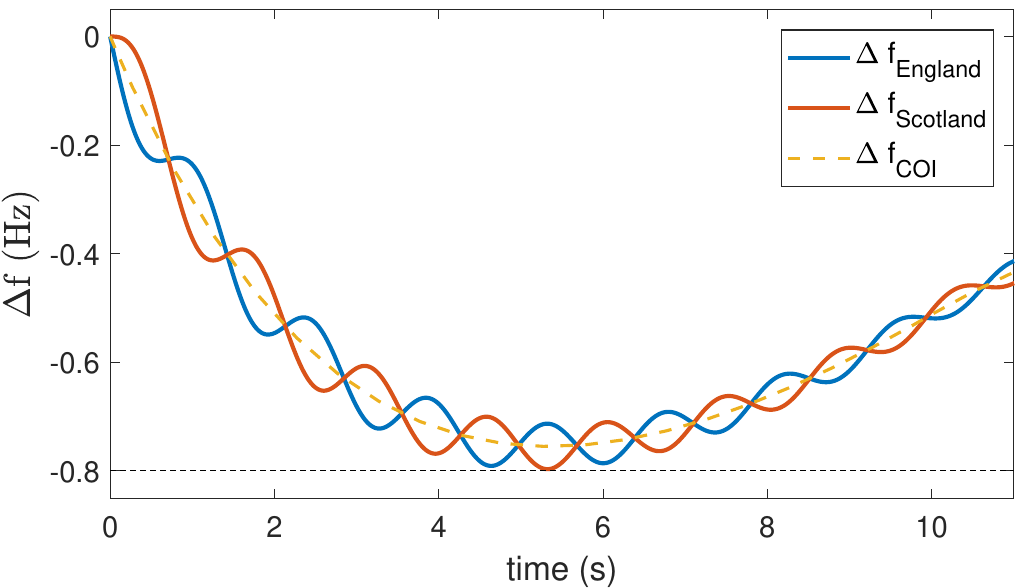}
        \caption{\footnotesize $\textrm{X}_{i,j}=40\Omega$}
        \vspace{2mm}
        \label{Small_impedance}
    \end{subfigure}%
    ~ 
    \newline 
    \begin{subfigure}[b]{0.5\textwidth}
        \centering
        \hspace*{-8mm}
        \includegraphics[width=3.1in]{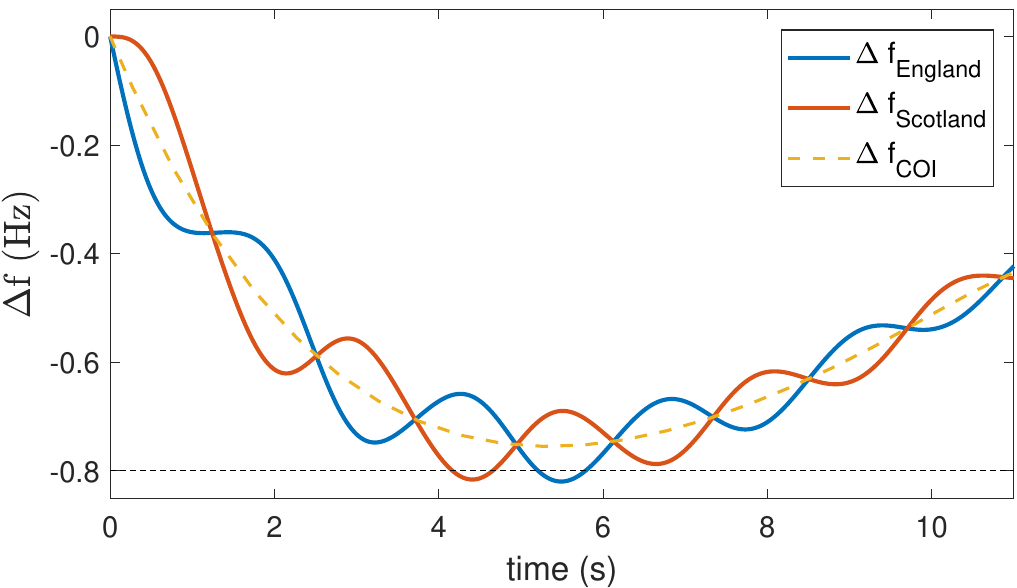}
        \caption{\footnotesize $\textrm{X}_{i,j}=100\Omega$}
        \label{Large_impedance}
    \end{subfigure}
    \caption{Impact of impedance of the transmission connection on the distinct regional frequencies. These cases consider a uniform inertia distribution, with both England and Scotland having the same inertia.}
    \label{Faults_differentImpedance}
\end{figure}

Both cases displayed in Fig.~\ref{Faults_differentImpedance}, which consider an even split of inertia, load damping and frequency response among England and Scotland, show a roughly equal COI frequency. However, some inter-area oscillations can be observed in both Figs. \ref{Small_impedance} and \ref{Large_impedance}, while the amplitude of these oscillations is notably higher in the case in Fig.~\ref{Large_impedance}. These oscillations clearly increase the absolute value of the RoCoF and further deteriorate the frequency nadir, which drops below $-0.8\textrm{Hz}$ (stability limit in GB) in Fig.~\ref{Large_impedance}.

High-impedance transmission lines are closely related to renewables, since the best renewable sources, particularly wind, are typically located in remote areas far from load centres. This fact implies that the electrical connection between these renewable generators and the rest of the grid will necessarily be weak given its great length (unless reinforcement of the network through, for example, parallel circuits is chosen, but this approach involves high investment costs).

\begin{figure}[!t]
	\centering
    \begin{subfigure}[b]{0.5\textwidth}
        \hspace*{-8mm}
        \centering
        \includegraphics[width=3.1in]{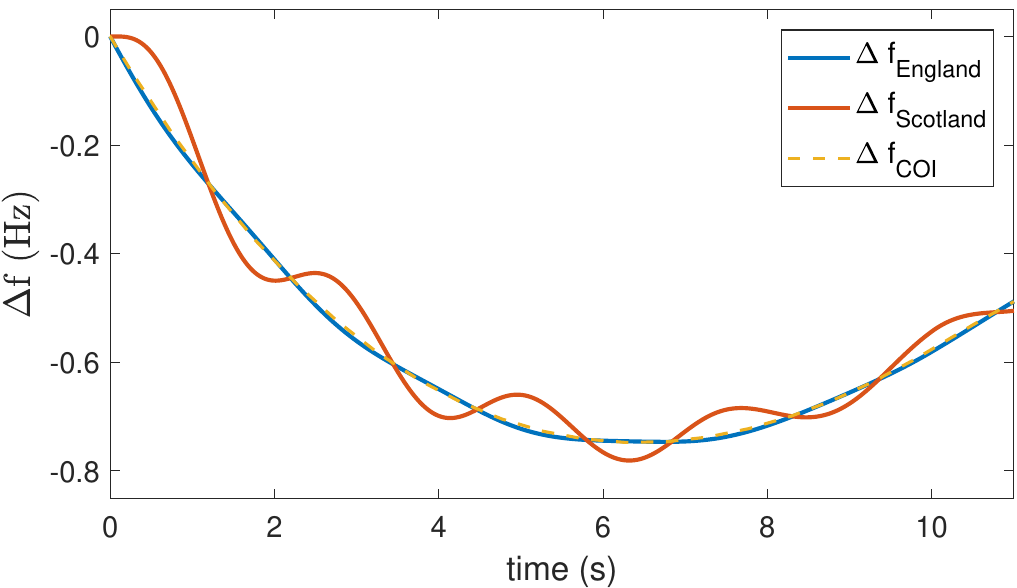}
        \caption{\footnotesize 1.8GW generation loss in England}
        \vspace{2mm}
        \label{Fault_England}
    \end{subfigure}%
    ~ 
    \newline 
    \begin{subfigure}[b]{0.5\textwidth}
        \hspace*{-8mm}
        \centering
        \includegraphics[width=3.1in]{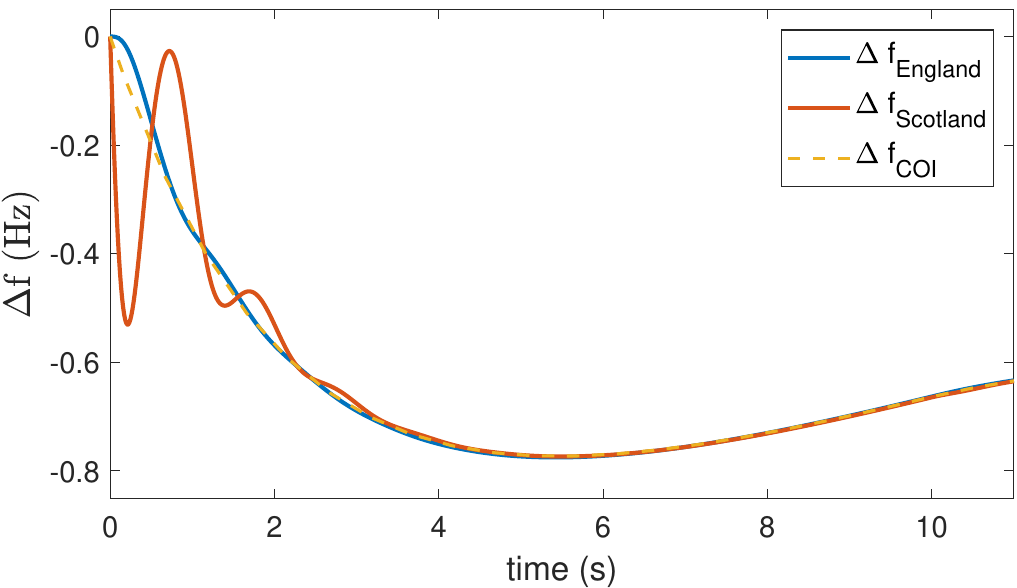}
        \caption{\footnotesize 0.6GW generation loss in Scotland}
        \label{Fault_Scotland}
    \end{subfigure}
    \caption{Impact of fault location on the distinct regional frequencies. System inertia is unevenly distributed in these cases, with 90\% located in England and 10\% in Scotland.}
    \label{Faults_EnglandScotland}
\end{figure}

These renewables located in remote regions also drive a non-uniform inertia distribution, as inertia in the regions dominated by renewable generation can be very low. The examples in Fig.~\ref{Faults_EnglandScotland} show that regional frequencies become significant when there are gradients of inertia in the system, particularly so when the generation outage occurs in the low-inertia region: the 0.6GW loss in Scotland displayed in Fig.~\ref{Fault_Scotland} exhibits extreme RoCoF in that region during the first few seconds, while a 1.8GW loss in England shows a similar behaviour to the COI model. The observation of highest oscillations in the first instants after the fault, for a fault happening in the low-inertia region was also reported in \cite{InmaMultiArea,DouglasWilsonSmartFreq}.

Oscillatory behaviour in different parts of an  electricity network is a complex phenomenon, not only influenced by the impedance of the transmission lines and the spatial distribution of system inertia, but where many factors come into play. The damping ratio of the electromechanical oscillations is also influenced by the excitation system of generators and loading level of lines. The results in this section show the importance of different regional frequencies, but no further assumption on the strength of the electrical connection or the inertia distribution is made in the remaining of this paper when deducing frequency-stability constraints. The impedance of the interconnection can be chosen in our model to match that of any particular system under study, while the location of inertia will be considered as a decision variable in our constraints.

\section{Dynamics of post-fault regional frequencies} \label{sectionDynamicModelMultiRegion}

This section describes the mathematical model of post-fault frequency evolution in multi-region systems, from which the conditions for regional frequency security will be obtained. For now let's consider the simplest case of 2 regions, as the generalisation for $N$ regions follows the same logic. The post-fault frequency in each region is described by this set of differential equations:
\begin{subequations} \label{2area_dynamic}
\begin{empheq}[left = \empheqlbrace]{align}
  & \; 2H_1\cdot\frac{\textrm{d}\Delta f_1(t)}{\textrm{d}t}+\textrm{D}_1\cdot \textrm{P}_1^\textrm{D}\cdot\Delta f_1(t)= \nonumber
  \\ & \qquad\qquad\qquad\qquad\quad \mbox{PFR}_1(t)-P_1^\textrm{L}+\Delta\textrm{P}^{\textrm{import}}_1(t) \tag{\ref{2area_dynamic}.1} \label{2area_dynamic_region1} \\[1pt]
  & \; 2H_2\cdot\frac{\textrm{d} \Delta f_2(t)}{\textrm{d}t}+\textrm{D}_2\cdot \textrm{P}_2^\textrm{D}\cdot\Delta f_2(t)= \nonumber
  \\ & \qquad\qquad\qquad\qquad\quad \mbox{PFR}_2(t)-P_2^\textrm{L}+\Delta\textrm{P}^{\textrm{import}}_2(t) \tag{\ref{2area_dynamic}.2} \label{2area_dynamic_region2}
\end{empheq}
\end{subequations}

Eqs. (\ref{2area_dynamic}) can be described as two coupled swing equations, in which the coupling term is the ac power exchange between regions. The coupling term $\Delta\textrm{P}^{\textrm{import}}_i(t)$ represents the deviation from the steady state power being imported to region $i$. That is, before the outage there was a certain amount of power being imported to region $i$; after the outage, the amount of power imported changes following the laws of energy conservation, so that a frequency equilibrium is restored in the whole network. $\Delta\textrm{P}^{\textrm{import}}_i(t)$ is the difference between the power imported after-outage and before-outage.

To obtain the mathematical description of the term $\Delta\textrm{P}^{\textrm{import}}_i(t)$, let's first consider the steady-state power transfer between two neighbouring buses, $i$ and $j$. This steady-state power transfer is described by the following equation \cite{KundurBook}:
\begin{equation}
\textrm{P}^\textrm{transfer}_{i,j}=\frac{\textrm{V}_i\textrm{V}_j}{\textrm{X}_{i,j}}\cdot\sin(\delta_i-\delta_j)
\end{equation}

\noindent Note that the line resistance is neglected as we consider $\textrm{X}_{i,j}\gg\textrm{R}_\textrm{i,j}$ (a typical assumption in transmission lines),
and that $\delta_i$, $\delta_j$ are bus voltage phase-angles referred to the rotating frame, which rotates at synchronous speed (50Hz in Europe). Therefore, these phase angles $\delta_i$, $\delta_j$ are time-invariant in steady state. The phase angles will however change in the event of a generation outage in the system, as the rotating speed of the generators deviates from 50Hz while they release the kinetic energy stored in their rotating masses.

The deviation from the steady state power transfer between buses $i$, $j$ after a fault can be calculated as:
\begin{multline}
\Delta \textrm{P}^\textrm{transfer}_{i,j}(t)  = \textrm{P}^\textrm{transfer}_{i,j}(t)-\textrm{P}^\textrm{transfer, ss}_{i,j} = \\ \frac{\textrm{V}_i\textrm{V}_j}{\textrm{X}_{i,j}}\cdot\sin\Big(\delta_i(t)-\delta_j(t)\Big) - \frac{\textrm{V}_i\textrm{V}_j}{\textrm{X}_{i,j}}\cdot\sin(\delta_i^\textrm{ss}-\delta_j^\textrm{ss})
\end{multline}

\noindent Linearising around the operating point $\delta_i=\delta^\textrm{ss}_i$, $\delta_j=\delta^\textrm{ss}_j$  \cite{KundurBook}:
\begin{equation}
\Delta \textrm{P}^\textrm{transfer}_{i,j}(t) = \frac{\textrm{V}_i\textrm{V}_j}{\textrm{X}_{i,j}}\cdot\cos(\delta_i^\textrm{ss}-\delta_j^\textrm{ss})\cdot\big[\Delta\delta_i(t)-\Delta\delta_j(t)\big]
\end{equation}

\noindent where $\Delta\delta_i(t)=\delta_i(t)-\delta_i^\textrm{ss}$.
We have assumed that the average voltage variation in the grid is small after the generation outage, an assumption also made in the literature \cite{RasoulPLossEstimate,SpatioTemporalMultiArea}.

Now, in order to express the power imported by any region in terms of electric frequency $\Delta f(t)$ rather than phase angle $\Delta\delta(t)$, the relation between these two magnitudes is given by:
\begin{equation}
\frac{\textrm{d}\delta(t)}{\textrm{d}t}=\frac{\textrm{d}\theta_m(t)}{\textrm{d}t}-\omega_s=\omega_m(t)-\omega_s=\Delta\omega(t)={2\pi}\Delta f(t)
\end{equation}

\noindent Where $\delta$ is assumed to be expressed in radians and $\Delta f$ in Hz. $\theta_m$ is the absolute angle of the synchronous generator, not referred to the rotating reference as $\delta$ is; $\omega_s$ is the synchronous rotational speed of the generators (equivalent to the 50Hz electrical frequency). Note that this relationship still holds when considering phase-angle deviations from steady state, $\Delta\delta(t)$:
\begin{equation}
\frac{\textrm{d}\Delta\delta(t)}{\textrm{d}t}=\frac{\textrm{d}\delta(t)}{\textrm{d}t}-\cancelto{0}{\frac{\textrm{d}\delta^\textrm{ss}}{\textrm{d}t}}=\Delta\omega(t)={2\pi}\Delta f(t)
\end{equation}

So the post-fault power imported by any region $i$ can be expressed as:
\begin{equation} \label{TransferRegions}
\Delta\textrm{P}^{\textrm{import}}_i(t)=-\sum_{j \in \mathcal{J}_i}{\textrm{T}_{i,j}}\left[\int_{0}^{t}\Delta f_i(\tau)d\tau-\int_{0}^{t}\Delta f_j(\tau)d\tau\right]
\end{equation}

Where $\mathcal{J}_i$ is the set of neighbour areas of area $i$. Note that a negative value of $\Delta\textrm{P}^{\textrm{import}}_i(t)$ means that region $i$ is exporting power to other regions. To make (\ref{TransferRegions}) and following expressions clearer, we have defined the electrical stiffness of the transmission line as \cite{ExpositoBook}:
\begin{equation}
    \textrm{T}_{i,j} = {2\pi}\cdot\frac{\textrm{V}_i\textrm{V}_j}{\textrm{X}_{i,j}}\cdot \cos(\delta_i^\textrm{ss}-\delta_j^\textrm{ss})
\end{equation}

The generalization for $N$ regions of the coupled swing equations, including the inter-region power transfer term, follows immediately:

\begin{subequations} \label{N_area_dynamics}
\begin{empheq}[left = \empheqlbrace]{align}
 & \; 2H_1\cdot\frac{\textrm{d} \Delta f_1(t)}{\textrm{d}t}+\textrm{D}_1\cdot \textrm{P}_1^\textrm{D}\cdot\Delta f_1(t)= \nonumber
  \\ & \qquad\qquad\qquad\qquad \textrm{PFR}_1(t)-P_1^\textrm{L}+\Delta\textrm{P}^{\textrm{import}}_1(t) \tag{\ref{N_area_dynamics}.1} \\[1pt]
 & \; 2H_2\cdot\frac{\textrm{d} \Delta f_2(t)}{\textrm{d}t}+\textrm{D}_2\cdot \textrm{P}_2^\textrm{D}\cdot\Delta f_2(t)= \nonumber
  \\ & \qquad\qquad\qquad\qquad \textrm{PFR}_2(t)-P_2^\textrm{L}+\Delta\textrm{P}^{\textrm{import}}_2(t)\tag{\ref{N_area_dynamics}.2}\\[-1pt]
 & \hspace{13.5em}\mathrel{\makebox[\widthof{=}]{\vdots}} \nonumber \\[-1pt]
 & \; 2H_n\cdot\frac{\textrm{d} \Delta f_n(t)}{\textrm{d}t}+\textrm{D}_n\cdot \textrm{P}_n^\textrm{D}\cdot\Delta f_n(t)= \nonumber
  \\ & \qquad\qquad\qquad\qquad \textrm{PFR}_n(t)-P_n^\textrm{L}+\Delta\textrm{P}^{\textrm{import}}_n(t) \tag{\ref{N_area_dynamics}.n}
\end{empheq}
\end{subequations}

\noindent where $\Delta\textrm{P}^{\textrm{import}}_i(t)$ is modelled as in (\ref{TransferRegions}).

The set of integro-differential equations (\ref{N_area_dynamics}) is the extension of the single-bus swing equation to consider distinct regions,
and therefore allows to consider geographical gradients in frequency within a power system. The power system is modelled as a set of regions that do not swing coherently but are connected through synchronous ac transmission lines. This model generalizes the uniform frequency model, as (\ref{N_area_dynamics}) reduces to the single-bus swing equation if the impedance of the transmission lines tends to zero. 

The groups of synchronous machines that exhibit a similar frequency behaviour can be identified through `coherency analysis’. The details of this technique can be found in reference \cite{CoherencyAnalysis}. Once these groups have been identified, this will determine the number of equations needed to build the system of eqs.~(\ref{N_area_dynamics}). Note that the `regions’ used in this paper are defined from the electromechanical coupling in a system, not necessarily corresponding to the `areas’ used for load frequency control. Therefore, intra-area oscillations can be considered in this formulation: a group of synchronous generators that swing against the rest of the system could be defined as a new region to be included in eqs.~(\ref{N_area_dynamics}).

\subsection{Solving the dynamic model: two-region case} \label{SectionSolving2area}

This section discusses how to obtain a solution for $\Delta f_1(t)$ and $\Delta f_2(t)$ as defined by (\ref{2area_dynamic}). Without loss of generality and for the sake of simplicity, we consider a single response service, Primary Frequency Response (PFR), modelled as:
\begin{subequations} \label{PFRmultiArea}
\begin{empheq}[left ={\textrm{PFR}_i(t)=\empheqlbrace}]{alignat=3}
 & \;\; \frac{R_i}{\textrm{T}_\textrm{g}} \cdot t \quad && \quad \mbox{if $t\leq\textrm{T}_\textrm{g}$} \\
 & \;\; R_i \quad && \quad \mbox{if $t> \textrm{T}_\textrm{g}$}
\end{empheq}
\end{subequations}

In order to obtain an analytical solution for the set of integro-differential equations (\ref{2area_dynamic}), we propose to first apply the Laplace transform so that a set of algebraic equations is obtained. Then, this set of algebraic equations can be solved analytically, and the inverse Laplace transform can be applied to obtain the solution in time domain. 

Applying the Laplace transform to (\ref{2area_dynamic}) we obtain:

\begin{subequations} \label{2area_Laplace}
\vspace*{-3mm}
\begin{empheq}[left = \empheqlbrace]{align}
  & \; 2H_1\cdot s\Delta F_1(s) +\textrm{D}_1\cdot \textrm{P}_1^\textrm{D}\cdot\Delta F_1(s)= \nonumber
  \\ & \qquad\qquad\qquad\qquad \mbox{PFR}_1(s)-\frac{P_1^\textrm{L}}{s}+\Delta\textrm{P}^{\textrm{import}}_1(s) \tag{\ref{2area_Laplace}.1}\\[1pt]
  & \; 2H_2 \cdot s\Delta F_2(s) +\textrm{D}_2\cdot \textrm{P}_2^\textrm{D}\cdot\Delta F_2(s)= \nonumber
  \\ & \qquad\qquad\qquad\qquad \mbox{PFR}_2(s)-\frac{P_2^\textrm{L}}{s}+\Delta\textrm{P}^{\textrm{import}}_2(s) \tag{\ref{2area_Laplace}.2}
\end{empheq}
\end{subequations}

\noindent where $\Delta F_i(s)=\Laplace\{\Delta f_i(t)\}$ and $\Delta\textrm{P}^{\textrm{import}}_i(s)$ is given by:
\begin{equation} \label{LaplaceImport}
\Delta\textrm{P}^{\textrm{import}}_i(s) =  -\sum_{j \in \mathcal{J}_i}\textrm{T}_{i,j}\frac{\Delta F_i(s)-\Delta F_j(s)}{s}
\end{equation}

\noindent Considering that PFR is still ramping up, which is the period of interest for RoCoF and nadir purposes:
\begin{equation}
\mbox{PFR}_i(s) = \frac{R_i}{\textrm{T}_\textrm{g}}\cdot\frac{1}{s^2}
\end{equation}


\addtocounter{equation}{4}

\begin{figure*}[!bh]
\vspace{-5pt}
\noindent\rule{\textwidth}{0.7pt}
    \begin{equation} \label{FirstSolutionLaplace}
    \Delta F_1(s) = \frac{- 2 H_2 P^\textrm{L}_1 \cdot s^3 + (\frac{2H_2R_1}{\textrm{T}_\textrm{g}} - \textrm{D}_2^{'}P^\textrm{L}_1)s^2 + (\frac{\textrm{D}_2^{'}R_1}{\textrm{T}_\textrm{g}} - P^\textrm{L}_1 \cdot \textrm{T}_{1,2})s + \frac{R \cdot \textrm{T}_{1,2}}{\textrm{T}_\textrm{g}}}{4H_1H_2 \cdot s^5 + 2(\textrm{D}_1^{'}H_2+\textrm{D}_2^{'}H_1)s^4 + (2H\cdot\textrm{T}_{1,2}+\textrm{D}_1^{'}\textrm{D}_2^{'})s^3 + \textrm{D}^{'}\textrm{T}_{1,2} \cdot s^2}
    \end{equation}
\vspace*{-15pt}
\end{figure*}

\begin{figure*}[!bp] 
    \begin{multline} \label{Laplace_partfrac}
    \Delta F_1(s) = \frac{R}{\textrm{D}^{'}\textrm{T}_\textrm{g}\cdot s^2} - \frac{2HR + \textrm{D}^{'}\textrm{T}_\textrm{g}P^\textrm{L}_1 - \frac{\textrm{D}_2^{'}(R_1\textrm{D}_2^{'}-R_2\textrm{D}_1^{'})}{\textrm{T}_{1,2}}}{(\textrm{D}^{'})^2\textrm{T}_\textrm{g}\cdot s}\\[-1pt]
    + \frac{\textrm{C}_1''s^2+\textrm{C}_2''s+\textrm{C}_3''}{\textrm{T}_\textrm{g}\textrm{T}_{1,2}(\textrm{D}^{'})^2 [4H_1H_2\cdot s^3 + 2(\textrm{D}_1^{'}H_2+\textrm{D}_2^{'}H_1)\cdot s^2 + (2H\cdot \textrm{T}_{1,2}+\textrm{D}_1^{'}\textrm{D}_2^{'})\cdot s + \textrm{D}^{'}\textrm{T}_{1,2}]}
    \end{multline}
\end{figure*}

\addtocounter{equation}{-6}

The set of algebraic equations (\ref{2area_Laplace}) where the unknowns are $\Delta F_i(s)$ can be solved analytically (e.g.~using a computer algebra system such as MATLAB or Maple), which gives the following solution:
\begin{equation} \label{GenericSolutionLaplace}
\Delta F_i(s) = \frac{\textrm{C}_1 \cdot s^3 + \textrm{C}_2 \cdot s^2 + \textrm{C}_3 \cdot s + \textrm{C}_4}{\textrm{C}_5 \cdot s^5 + \textrm{C}_6 \cdot s^4 + \textrm{C}_7 \cdot s^3 + \textrm{C}_8 \cdot s^2}
\end{equation}

\noindent Where the constants $\textrm{C}_j$ are all real numbers, which are functions of the system parameters that appear in (\ref{2area_Laplace}) such as $H_1$, $H_2$, $P^\textrm{L}_1$, etc. These functions of system parameters are significantly convoluted as shown later in this section, hence the use of these generic constants $\textrm{C}_j$ to represent them in a simple way.
 
In order to apply the inverse Laplace transform and then obtain the time-domain solution $\Delta f_i(t)$, (\ref{GenericSolutionLaplace}) must be decomposed in partial fractions:
\begin{multline} \label{FirstPartFrac}
\Delta F_i(s) = \frac{\textrm{C}_1'}{s^2} + \frac{\textrm{C}_2'}{s} + \frac{\textrm{C}_3'\cdot s^2+\textrm{C}_4'\cdot s+1}{\textrm{C}_5'\cdot s^3+\textrm{C}_6'\cdot s^2+\textrm{C}_7'\cdot s+1} = \\[0pt] \frac{\textrm{C}_1'}{s^2} + \frac{\textrm{C}_2'}{s} + \sum_{k=1}^{3}\frac{\textrm{Z}_k}{s-\textrm{z}_k}
\end{multline}

\noindent Note that these constants $\textrm{C}_j'$ are again generic functions of the system parameters but different to the $\textrm{C}_j$ used in (\ref{GenericSolutionLaplace}). Terms $\textrm{Z}_k$ and $\textrm{z}_k$ are generic complex numbers. The factorisation performed in (\ref{FirstPartFrac}) makes use of the \textbf{Fundamental Theorem of Algebra}, which states that every polynomial of degree $n$ has $n$ roots. The inverse Laplace transform of (\ref{FirstPartFrac}) gives the time-domain solution:
\vspace{-0.5mm}
\begin{multline} \label{FirstTimeSolution}
\hspace*{-3.5mm} \Delta f_i(t) = \Laplace^{-1}\left\{\Delta F_i(s) \right\} = \Laplace^{-1}\left\{\frac{\textrm{C}_1'}{s^2} + \frac{\textrm{C}_2'}{s} + \sum_{k=1}^{3}\frac{\textrm{Z}_k}{s-\textrm{z}_k}\right\} \\[0pt]
= \textrm{C}_1' \cdot t + \textrm{C}_2' + \sum_{k=1}^{3} \textrm{Z}_k \cdot \textrm{e}^{\textrm{z}_k t}
\end{multline}

To fully understand (\ref{FirstTimeSolution}), let's consider the two possibilities for the roots of the last denominator in (\ref{FirstPartFrac}), $\textrm{z}_k$, : 1) All $\textrm{z}_k$ are real; or 2) $\textrm{z}_1$ is real while $\textrm{z}_2$ and $\textrm{z}_3$ are complex conjugates. This conclusion comes from the \textbf{Complex Conjugate Root Theorem}, which states that if a polynomial in one variable with real coefficients (as is the case in (\ref{FirstPartFrac})) has a complex root, then the conjugate is also a root. 
\begin{prop} \label{Proposition_2_areas}
For the range of values in a realistic power system, the third-order denominator in (\ref{FirstPartFrac}) has one real root and two complex conjugate roots, i.e.~$\textrm{z}_1$ is real while $\textrm{z}_2$ and $\textrm{z}_3$ are complex conjugates in (\ref{FirstPartFrac}) and (\ref{FirstTimeSolution}). 
\end{prop}
The validity of this proposition is supported by the evidence reported in the literature \cite{InmaMultiArea,DouglasWilsonSmartFreq,OverbyeLocationalInertia,GoranAndersson1} and in the dynamic simulations in Section \ref{SectionDynamicSimulations}, which show that the post-fault regional frequencies behave as the COI plus certain inter-area oscillations. As is shown in the next paragraphs, such dynamic behaviour corresponds to one real root and two complex conjugate roots for the denominator in (\ref{FirstPartFrac}). Note that the alternative to Proposition~\ref{Proposition_2_areas} would be that the third-order denominator has three real roots, which would give rise to a sum of three exponential functions in time-domain (with no oscillations); this behaviour would not be consistent with that of realistic power systems.

Considering Proposition~\ref{Proposition_2_areas}, the last term in (\ref{FirstTimeSolution}) can be expressed as:
\vspace{-1mm}
\begin{multline} \label{DecompositionComplexRoots}
\sum_{k=1}^{3} \textrm{Z}_k \cdot \textrm{e}^{\textrm{z}_k t} = \textrm{Z}_1 \cdot \textrm{e}^{\textrm{z}_1 t} + \textrm{Z}_2 \cdot \textrm{e}^{\textrm{z}_2 t} + \overline{\textrm{Z}_2} \cdot \textrm{e}^{\overline{\textrm{z}_2} t} \\
= \textrm{Z}_1 \cdot \textrm{e}^{\textrm{z}_1 t} + (\textrm{a}_2+j\textrm{b}_2) \cdot \textrm{e}^{(\textrm{a}_2+j\textrm{b}_2) t} + (\textrm{a}_2-j\textrm{b}_2) \cdot \textrm{e}^{(\textrm{a}_2-j\textrm{b}_2) t} \\
= \textrm{Z}_1 \cdot \textrm{e}^{\textrm{z}_1 t} + \textrm{e}^{\textrm{a}_2 t} \left[ \textrm{a}_2 \left(\textrm{e}^{j\textrm{b}_2 t}+\textrm{e}^{-j\textrm{b}_2 t} \right) + j\textrm{b}_2 \left(\textrm{e}^{j\textrm{b}_2 t}-\textrm{e}^{-j\textrm{b}_2 t} \right) \right] \\
= \textrm{Z}_1 \cdot \textrm{e}^{\textrm{z}_1 t} + \textrm{e}^{\textrm{a}_2 t} \left[\textrm{a}_2 \cdot 2 \cos(\textrm{b}_2 t) - \textrm{b}_2 \cdot 2 \sin(\textrm{b}_2 t) \right]
\end{multline}

\noindent Where $\textrm{a}_2$ and $\textrm{b}_2$ are generic real numbers. Note that (\ref{DecompositionComplexRoots}) shows that the last term in (\ref{FirstTimeSolution}) corresponds to an exponential function plus certain oscillations in time-domain.

Since a third-order polynomial can be factorised algebraically (for example, using Cardano's formula), we now describe the procedure to do so. First, let's consider the explicit values for the generic constants $\textrm{C}_j$ in (\ref{GenericSolutionLaplace}) and $\textrm{C}_j'$ in (\ref{FirstPartFrac}) for region 1 (the expression for region 2 is omitted but equivalent): see eq. (\ref{FirstSolutionLaplace}).

\widowpenalty=0 
Eq. (\ref{FirstSolutionLaplace}) has been simplified using: $H=H_1+H_2$, $R=R_1+R_2$, $\textrm{D}_i^{'}=\textrm{D}_i\cdot\textrm{P}^\textrm{D}_i$ and $\textrm{D}{'}=\textrm{D}_1\cdot\textrm{P}^\textrm{D}_2+\textrm{D}_1\cdot\textrm{P}^\textrm{D}_2$. 
Before applying the inverse Laplace transform and then obtain the time-domain solution $\Delta f_1(t)$, (\ref{FirstSolutionLaplace}) must be decomposed in partial fractions as was done to obtain (\ref{FirstPartFrac}) from (\ref{GenericSolutionLaplace}). An algebraic partial-fraction decomposition can be done using a computer algebra system, giving eq.~(\ref{Laplace_partfrac}).
The value of the terms `C' in eq.~(\ref{Laplace_partfrac}) is:

\addtocounter{equation}{2}

\begin{equation}
\textrm{C}_1'' = 4H_1H_2(P^\textrm{L}_1\textrm{T}_\textrm{g}\textrm{T}_{1,2}\textrm{D}^{'} + 2HR\cdot\textrm{T}_{1,2} - R_1(\textrm{D}_2^{'})^2 + R_2\textrm{D}_1^{'}\textrm{D}_2^{'})
\end{equation}

\vspace*{-4mm}
\begin{multline}
\hspace{-3.5mm}\textrm{C}_2'' = 2[P^\textrm{L}_1\textrm{T}_\textrm{g}\textrm{T}_{1,2}\textrm{D}_2^{'}(H_1-H_2)(\textrm{D}_1^{'}+\textrm{D}_2^{'}) \\ 
+ 2R\cdot\textrm{T}_{1,2}(\textrm{D}_2^{'}H_1^2+\textrm{D}_1^{'}H_2^2)
+ \textrm{D}_1^{'}(\textrm{D}_2^{'})^2(H_1R_2-H_2R_1) \\ 
+ (\textrm{D}_1^{'})^2\textrm{D}_2^{'}H_2R_2
- (\textrm{D}_2^{'})^3H_1R_1]
\end{multline}

\vspace*{-4mm}
\begin{multline} \label{lastConstant2area}
\textrm{C}_3'' = 4 H^2 R (\textrm{T}_{1,2})^2 
+2H\cdot\textrm{D}^{'}P^\textrm{L}_1\textrm{T}_\textrm{g}(\textrm{T}_{1,2})^2 \\
- (\textrm{D}_2^{'})^2\textrm{T}_{1,2}[P^\textrm{L}_1\textrm{T}_\textrm{g}\textrm{D}^{'}+2H_1(2R_1+R_2)] \\
+2\textrm{D}_1^{'}\textrm{D}_2^{'}\textrm{T}_{1,2}(H_1R_2+2H_2R_1+H_2R_2) \\
-2(\textrm{D}_1^{'})^2\textrm{T}_{1,2}H_2R_2
+\textrm{D}_1^{'}(\textrm{D}_2^{'})^2(R_2\textrm{D}_1^{'}-R_1\textrm{D}_2^{'})
\end{multline}

At this point, it is revealing to compare (\ref{Laplace_partfrac}) with the partial-fraction decomposition of the uniform frequency model, which was deduced in \cite{FeiStochastic}:

\begin{equation} \label{LaplaceCOI}
\Delta F_\textrm{COI}(s) = \frac{R}{\textrm{D}^{'}\textrm{T}_\textrm{g}\cdot s^2} + \frac{2HR + \textrm{D}^{'}\textrm{T}_\textrm{g}P^\textrm{L}}{(\textrm{D}^{'})^2\textrm{T}_\textrm{g}}\left( \frac{1}{s + \frac{\textrm{D}^{'}}{2H}} - \frac{1}{s} \right)
\end{equation}

By applying the inverse Laplace transform to (\ref{LaplaceCOI}), the $\frac{1}{s^2}$ term would give rise to a linear function in time-domain, the $\frac{1}{s}$ term to a constant in time-domain and the $\frac{1}{s+\frac{\textrm{D}'}{2H}}$ term to an exponential function in time-domain:

\begin{equation} \label{eq:COIfrequency}
\Delta f_\textrm{COI}(t) = \frac{R}{\textrm{D}^{'}\textrm{T}_\textrm{g}} \cdot t + \frac{2HR + \textrm{D}^{'}\textrm{T}_\textrm{g}P^\textrm{L}}{(\textrm{D}^{'})^2\textrm{T}_\textrm{g}} \left( \textrm{e}^{-\frac{\textrm{D}^{'}}{2H} \cdot t} - 1 \right)
\end{equation}


By comparing (\ref{Laplace_partfrac}) to (\ref{LaplaceCOI}), one can notice that the $\frac{1}{s^2}$ term is the same in both expressions, and the $\frac{1}{s}$ in (\ref{Laplace_partfrac}) reduces to the one in (\ref{LaplaceCOI}) as $\textrm{T}_{1,2}$ tends to infinity, i.e.~the impedance of the interconnecting line tends to zero and therefore the two regions become effectively one. The comparison of the $\frac{1}{s+\frac{\textrm{D}'}{2H}}$ term in (\ref{LaplaceCOI}) with the last fraction in (\ref{Laplace_partfrac}) is not as straightforward, but the evidence reported in the literature and in the simulations conducted in Section \ref{SectionDynamicSimulations} shed some light for the comparison: the post-fault frequency in a multi-region system exhibits the behaviour of the COI plus certain inter-area oscillations. This observation is consistent with (\ref{DecompositionComplexRoots}), which demonstrates that the last term in (\ref{Laplace_partfrac}) corresponds to an exponential function plus oscillations in time-domain. Therefore the exponential time-domain term rising from the $\frac{1}{s+\frac{\textrm{D}'}{2H}}$ term in (\ref{LaplaceCOI}) is roughly equal to the exponential term in (\ref{DecompositionComplexRoots}).

In conclusion, the deductions in this section have shown that for practical purposes in any realistic power system, the post-fault frequency evolution in a two-region system is of the form:
\begin{equation} \label{SolutionFreq_2areas}
\Delta f_i(t) \approx \Delta f_{\textrm{COI}}(t) + \textrm{e}^{-\textrm{a} t} \textrm{A}_i \sin(\omega  t+\phi_i)+\textrm{C}_i
\end{equation}
Note that the approximation in (\ref{SolutionFreq_2areas}) lies in the COI term, as discussed in the paragraph above.

The deductions in this section have allowed to obtain the \textit{mathematical structure} of the post-fault frequency evolution in a two-region system, given by eq.~(\ref{SolutionFreq_2areas}). This expression will be used as an initial step to obtain conditions that guarantee frequency security in a two-region system, in Section \ref{Conditions2Region}.

Finally, it is worth noting that since the third-order denominator in (\ref{Laplace_partfrac}) can be factorised analytically, in principle the exact expression that has been approximated in (\ref{SolutionFreq_2areas}) could be obtained analytically. However, this is not done in this paper for two reasons: 1) given that the expressions (\ref{Laplace_partfrac}) through (\ref{lastConstant2area}) are cumbersome, the solution would include highly non-linear functions which would make frequency-security constraints very inefficient for being implemented in optimisation routines (if even possible to obtain such constraints); and 2) deducing this analytical solution would only be possible for a two-region system, and not for a system with three or more regions, as is proved in the next section.

\subsection{Solving the dynamic model: N-region case} \label{Freq_i_N_area}
This section discusses how to obtain a solution for every region's post-fault frequency, $\Delta f_i(t)$, in an $N$-region system as defined by (\ref{N_area_dynamics}). 
The Laplace transform of this set of integro-differential eqs.~(\ref{N_area_dynamics}) is given by:
\begin{subequations} \label{Narea_Laplace}
\begin{empheq}[left = \empheqlbrace]{align}
  & \; 2H_1\cdot s\Delta F_1(s) +\textrm{D}_1\cdot \textrm{P}_1^\textrm{D}\cdot\Delta F_1(s)= \nonumber
  \\ & \qquad\qquad\qquad\qquad \mbox{PFR}_1(s)-\frac{P_1^\textrm{L}}{s}+\Delta\textrm{P}^{\textrm{import}}_1(s) \tag{\ref{Narea_Laplace}.1}\\[5pt]
  & \; 2H_2 \cdot s\Delta F_2(s) +\textrm{D}_2\cdot \textrm{P}_2^\textrm{D}\cdot\Delta F_2(s)= \nonumber
  \\ & \qquad\qquad\qquad\qquad \mbox{PFR}_2(s)-\frac{P_2^\textrm{L}}{s}+\Delta\textrm{P}^{\textrm{import}}_2(s) \tag{\ref{Narea_Laplace}.2}\\[1pt]
 & \hspace{13.5em}\mathrel{\makebox[\widthof{=}]{\vdots}} \nonumber \\[1pt]
 & \; 2H_n\cdot s\Delta F_n(s)+\textrm{D}_n\cdot \textrm{P}_n^\textrm{D}\cdot\Delta F_n(s)= \nonumber
  \\ & \qquad\qquad\qquad\qquad \textrm{PFR}_n(s)-\frac{P_n^\textrm{L}}{s}+\Delta\textrm{P}^{\textrm{import}}_n(s) \tag{\ref{Narea_Laplace}.n}
\end{empheq}
\end{subequations}


\noindent where $\Delta\textrm{P}^{\textrm{import}}_i(s)$ is given by (\ref{LaplaceImport}). The solution for a generic region $i$ of this set of algebraic eqs. is:

\begin{equation} \label{GeneralFormulaNRegions}
\Delta F_i(s) = \frac{\textrm{C}_1}{s^2} - \frac{\textrm{C}_2}{s} + \frac{\textrm{polynomial in } s \textrm{ of degree } 2\cdot N-2}{\textrm{polynomial in } s \textrm{ of degree } 2\cdot N-1} 
\end{equation}

In order to apply the inverse Laplace transform to (\ref{GeneralFormulaNRegions}) and then obtain the time-domain solution $\Delta f_i(t)$, the last term in (\ref{GeneralFormulaNRegions}) must be decomposed in partial fractions. However, since it is not possible to algebraically factorise polynomials of degree higher than four as proved by the \textbf{Abel--Ruffini Theorem}, it is not possible to decompose in partial fractions the last term in (\ref{GeneralFormulaNRegions}). Therefore, it is not possible to obtain an analytical solution in time-domain for $\Delta f_i(t)$ for a general power system with $N$ regions, i.e.~a system with three or more regions. 
\begin{prop} \label{Proposition_N_areas}
For the range of values in a realistic power system, the denominator of the last term in (\ref{GeneralFormulaNRegions}) always has one real root and `$2N-2$' complex conjugate roots.  
\end{prop}
This proposition makes use of the Fundamental Theorem of Algebra and the Complex Conjugate Roots theorem, see the previous Section \ref{SectionSolving2area} discussing two-region systems for more details.
Proposition~\ref{Proposition_N_areas} is consistent with the studies reported in the literature \cite{InmaMultiArea,DouglasWilsonSmartFreq,OverbyeLocationalInertia,GoranAndersson1}, which show post-fault frequencies behaving as the COI plus certain inter-area oscillations. Note that the alternative to Proposition~\ref{Proposition_N_areas} would be a sum of real roots, giving rise to a sum of exponential functions in time-domain, which is not consistent with the behaviour of realistic power systems.

Therefore, although no analytical solution can be obtained for $\Delta f_i(t)$ for a general system with $N$ regions, Proposition~\ref{Proposition_N_areas} gives information about the \textit{mathematical structure} of $\Delta f_i(t)$, in a similar fashion as was deduced in Section \ref{SectionSolving2area} for a two-region system.

As discussed in Section \ref{SectionSolving2area}, the real root of the denominator in of the last term in (\ref{GeneralFormulaNRegions}) gives rise to an exponential in time-domain, which added to the first two terms in (\ref{GeneralFormulaNRegions}) corresponds to the post-fault frequency of the COI; the $2N-2$ complex conjugate roots give rise to $N-1$ oscillations in time-domain, as was discussed for the two complex conjugate roots of the two-region system considered in Section \ref{SectionSolving2area}. Therefore, the solution for the post-fault frequency evolution of any region $i$ in an $N$-region system is:

\begin{equation} \label{SolutionFreq_N_areas}
\Delta f_i(t) \approx \Delta f_{\textrm{COI}}(t) + \sum_{j=1}^{N-1} \textrm{e}^{-\textrm{a}_j t} \textrm{A}_j \sin(\omega_j t+\phi_j)+\textrm{C}_j
\end{equation}

\section{Conditions for frequency security: 2 regions} \label{Conditions2Region}

In this section we propose a method to obtain conditions for respecting RoCoF, nadir and quasi-steady-state (q-s-s) limits in any region $i$. To illustrate this method, let's start by considering the time-evolution of frequency deviation in any of the two regions, $\Delta f_i (t)$, given by eq.~(\ref{SolutionFreq_2areas}).
It is not possible to obtain analytical expressions for the coefficients `$\textrm{a}$', `$\textrm{A}_i$' and `$\omega$' as functions of the system parameters ($H_1$, $H_2$, $P^\textrm{L}$, etc.); however, using the information provided by (\ref{SolutionFreq_2areas}), a numerical approach can be used to estimate these coefficients. 

This approach consists on numerically solving the set of differential eqs.~(\ref{N_area_dynamics}) for several possible operating points of a power system (therefore considering only plausible values for $H_1$, $H_2$, $P^\textrm{L}$ and the rest of the system variables), and using a regression technique on the solution samples to estimate the coefficients `$\textrm{a}$', `$\textrm{A}_i$' and `$\omega$'. In summary, this mixed analytical-numerical approach goes as far as possible using analytical techniques and completes the task of obtaining frequency-security constraints by using numerical techniques to estimate the remaining parameters that cannot be obtained analytically. This method is driven by the two following goals: 1) to always guarantee frequency security, therefore becoming conservative if necessary (although it will be shown that the conservativeness introduced is not significant); and 2) to obtain linear constraints that would provide computational efficiency for the optimisation problem in which they would eventually be implemented.

\subsection{RoCoF constraint: analytical deduction} \label{SecRoCoF2region}

Taking into consideration (\ref{SolutionFreq_2areas}) the RoCoF in region $i$ is:
\begin{equation} \label{RoCoF_function_i}
\textrm{RoCoF}_i(t) = \textrm{RoCoF}_\textrm{COI}(t) + \textrm{RoCoF}_{\textrm{oscillations}_i}(t)
\end{equation}

\noindent where $\textrm{RoCoF}_{\textrm{oscillations}_i}(t)$ is the derivative of the attenuated oscillations in (\ref{SolutionFreq_2areas}):
\begin{equation}
\textrm{RoCoF}_{\textrm{oscillations}_i}(t)=\pm e^{-\textrm{a} t}\textrm{A}_i[\omega \cos(\omega t+\phi_i)-\textrm{a}\cdot\sin(\omega t +\phi_i)]
\end{equation}


\noindent Since this $\textrm{RoCoF}_{\textrm{oscillations}_i}(t)$ is a non-convex function, then $\textrm{RoCoF}_i(t)$ in eq.~(\ref{RoCoF_function_i}) is also a non-convex function and therefore it is not possible to obtain its global maximum analytically. We make use of the following inequality to come up with a conservative estimation of this maximum for $\textrm{RoCoF}_{i}(t)$:
\begin{multline}
\sup\{\textrm{RoCoF}_i(t)\mid t \geq 0 \} \; \leq \;  \sup\{\textrm{RoCoF}_\textrm{COI}(t)\mid t \geq 0 \} \\ + \sup\{\textrm{RoCoF}_{\textrm{oscillations}_i}(t)\mid t \geq 0 \}
\end{multline}


The $\sup\{\textrm{RoCoF}_{\textrm{oscillations}_i}(t)\mid t \geq 0 \}$ cannot be obtained analytically, because $\textrm{RoCoF}_{\textrm{oscillations}_i}(t)$ is a non-convex function as mentioned before, but it can be overestimated by:
\begin{equation} \label{Supremum_Rocof_osc}
\sup\{\textrm{RoCoF}_\textrm{oscillations}(t)\mid t \geq 0 \} \leq \textrm{A}_i\cdot \omega
\end{equation}


\noindent Expression (\ref{Supremum_Rocof_osc}) has been obtained by neglecting the attenuation term $e^{-\textrm{a} t}$, therefore assuming that the oscillations are not attenuated.

In conclusion, the constraints that guarantee the RoCoF to be within specified limits in a two-region system are:

\begin{subequations} \label{eq:RoCoF_analytical_TwoRegions}
\begin{empheq}{align}
  & \left|\mbox{RoCoF}_1\right| = \frac{P^\textrm{L}}{2 (H_1+H_2)}+\textrm{A}_1\cdot\omega \leq \mbox{RoCoF}_{\textrm{max}}  \\[0pt]
  & \left|\mbox{RoCoF}_2\right| = \frac{P^\textrm{L}}{2 (H_1+H_2)}+\textrm{A}_2\cdot\omega \leq \mbox{RoCoF}_{\textrm{max}}
\end{empheq}
\end{subequations}

\noindent where $\textrm{A}_i$ and $\omega$ are dependent on the operating condition of the system. 
They can be estimated using a linear function on simulation samples, which yields linear RoCoF constraints. This is discussed in detail in Part II of this paper.

\subsection{Nadir constraint: analytical deduction} \label{SecNadir2region}

In a similar fashion as for the RoCoF constraints discussed in the previous Section \ref{SecRoCoF2region}, it is not possible to deduce a purely analytical nadir constraint from the solution of $\Delta f_i$ given by (\ref{SolutionFreq_2areas}), as it would involve finding the global minimum of that non-convex function (\ref{SolutionFreq_2areas}) (see Fig.~\ref{Faults_differentImpedance} for graphical evidence of the existence of several local minima in the post-fault frequency evolution of a region). Deducing the mathematical expression for this global minimum would be necessary to obtain the expression for $t^*$,
which is needed to obtain an analytical nadir constraint $\left| \Delta f_{\textrm{nadir}} \right| = \left| \Delta f(t=t^*)  \right| \leq \Delta f_{\textrm{max}}$.

To overcome this limitation, we consider the energy equilibrium, instead of the power equilibrium, at the time of the frequency nadir.
Assuming that the generation outage occurs in region 2, the nadir constraint formulated as an energy equilibrium takes this form:
\begin{multline} \label{nadir_region2}
\underbrace{\vphantom{\frac{R}{\textrm{T}}}P^{\textrm{L}}_2\cdot\textrm{t}^*}_{\substack{\text{Energy} \\ \text{``lost"}}}
\hspace*{2mm} \leq \hspace*{2mm} 
\underbrace{\underbrace{\vphantom{\frac{R_1}{\textrm{T}_\textrm{g}}} 2H_2\Delta f_\textrm{max}}_{\substack{\text{Energy} \\ \text{contribution} \\ \text{from inertia}}} +
\underbrace{\strut \frac{R_2}{\textrm{T}_\textrm{g}}\frac{(\textrm{t}^*)^2}{2}}_{\substack{\text{Energy} \\ \text{contribution} \\ \text{from PFR}}} + 
\underbrace{\vphantom{\frac{R_1}{\textrm{T}_\textrm{g}}} \textrm{D}_2\textrm{P}^\textrm{D}_2 \int_{0}^{\textrm{t}^*} \Delta f_2(t) \ \textrm{d}t}_{\substack{\text{Energy contribution} \\ \text{from damping}}}
}_\text{Maximum admissible energy ``injected"} \\
+
\underbrace{\underbrace{\vphantom{\frac{R_1}{\textrm{T}_\textrm{g}}}\textrm{T}_{1,2} \int_{0}^{\textrm{t}^*}\int_{0}^{t}\left[\Delta f_1(\tau)-\Delta f_2(\tau)\right] \textrm{d}\tau \ \textrm{d}t}_{\substack{\text{Energy imported} \\ \text{from region 1}}}
}_\text{Maximum admissible energy ``injected"}
\end{multline}
This expression has been obtained by integrating the swing equation for region 2, i.e.~integrating (\ref{2area_dynamic_region2}). In plain words, constraint (\ref{nadir_region2}) enforces that \textit{the energy ``lost" must be lower than the maximum admissible energy ``injected"}. This maximum admissible energy ``injected" is limited by the maximum kinetic energy that can be extracted from the rotating masses in the system so that frequency does not drop below $\Delta f_\textrm{max}$. It is clear from (\ref{nadir_region2}) that if the energy related to the inertia term is higher, the frequency deviation would be higher than $\Delta f_\textrm{max}$. The energy contribution from PFR, damping and imports from the non-faulted region are also accounted for, and the sum of these terms cannot be lower than the energy ``lost" from the generation outage $P^\textrm{L}_2$.


For region 1, the nadir constraint obtained by integrating (\ref{2area_dynamic_region1}) is as follows:
\begin{multline} \label{nadir_region1}
\hspace*{-3mm}
\underbrace{\textrm{T}_{1,2} \int_{0}^{\textrm{t}^*}\int_{0}^{t}\left[\Delta f_1(\tau)-\Delta f_2(\tau)\right] \textrm{d}\tau \ \textrm{d}t}_\text{Energy ``lost"}
\hspace*{2mm} \leq \hspace*{2mm} \\[0pt]
\underbrace{\underbrace{\vphantom{\frac{R_1}{\textrm{T}_\textrm{g}}} 2H_1\Delta f_\textrm{max}}_{\substack{\text{Energy} \\ \text{contribution} \\ \text{from inertia}}} +
\underbrace{\frac{R_1}{\textrm{T}_\textrm{g}}\frac{(\textrm{t}^*)^2}{2}}_{\substack{\text{Energy} \\ \text{contribution} \\ \text{from PFR}}} + 
\underbrace{\vphantom{\frac{R_1}{\textrm{T}_\textrm{g}}} \textrm{D}_1\textrm{P}^\textrm{D}_1 \int_{0}^{\textrm{t}^*}\Delta f_1(t) \ \textrm{d}t}_{\substack{\text{Energy contribution} \\ \text{from damping}}}
}_\text{Maximum admissible energy ``injected"}
\end{multline}
Which again states that the energy ``lost" (in this case energy sent to the outaged region 2 through the transmissions corridors) cannot be higher than the maximum admissible energy that can be extracted from inertia (while considering the energy contribution from PFR, damping and imports).

The integrals that include the terms `$\Delta f_1(t)$' and `$\Delta f_2(t)$' must be estimated in the above constraints, for which we propose to use a numerical methodology. This approach allows to overcome the fact that no analytical solution exists for `$\Delta f_i(t)$' as a function of the system operating state ($H_1$, $H_2$, $P_\textrm{2}^L$, etc.), as demonstrated in Section~\ref{SectionSolving2area}. This numerical methodology, explained in detail in Part II of this paper, estimates these terms using a linear function of the system operating state so that the final form of constraints (\ref{nadir_region2}) and (\ref{nadir_region1}) is linear.

Finally, since it is not possible to obtain the analytical expression for $t^*$
(as explained before, it would involve finding the global minimum of a non-convex function), we propose to discretise (\ref{nadir_region2}) and (\ref{nadir_region1}) using several time-intervals and apply conditional constraints. This discretisation is based on the idea that the nadir has been reached if the power delivered by PFR and damping is higher than the generation loss (since the power delivered by inertia at the nadir is zero). For $n$ segments in the discretisation of time $t \in [t_1,t_2,\dots,t_{n-1},\textrm{T}_\textrm{g}]$, the nadir constraints for a two-region system would take this shape:

\vspace*{2mm}

if $\quad \frac{R}{\textrm{T}_\textrm{g}} \cdot t_1 > P^\textrm{L}_2 - \textrm{D}\cdot\textrm{P}_\textrm{D}\cdot \Delta f_\textrm{max} \quad$ then enforce:
\vspace{-3mm}
\begin{multline} 
\textrm{T}_{1,2} \int_{0}^{\textrm{t}_1}\int_{0}^{t}\left[\Delta f_1(\tau)-\Delta f_2(\tau)\right] \textrm{d}\tau \ \textrm{d}t
\hspace*{2mm} \leq \hspace*{2mm} \\[-3pt]
2H_1\Delta f_\textrm{max} +
\frac{R_1}{\textrm{T}_\textrm{g}}\frac{(\textrm{t}_1)^2}{2} + 
\textrm{D}_1\textrm{P}^\textrm{D}_1 \int_{0}^{\textrm{t}_1}\Delta f_1(t) \ \textrm{d}t
\end{multline}
\vspace{-6mm}
\begin{multline} 
P^\textrm{L}_2\cdot t_1 
\hspace*{2mm} \leq \hspace*{2mm}
2H_2\Delta f_\textrm{max} +
\frac{R_2}{\textrm{T}_\textrm{g}}\frac{(\textrm{t}_1)^2}{2} + 
\textrm{D}_2\textrm{P}^\textrm{D}_2 \int_{0}^{\textrm{t}_1}\Delta f_2(t) \ \textrm{d}t \\[-2pt]
+ \textrm{T}_{1,2} \int_{0}^{\textrm{t}_1}\int_{0}^{t}\left[\Delta f_1(\tau)-\Delta f_2(\tau)\right] \textrm{d}\tau \ \textrm{d}t
\end{multline}

else if $\quad \frac{R}{\textrm{T}_\textrm{g}} \cdot t_2 > P^\textrm{L}_2 - \textrm{D}\cdot\textrm{P}_\textrm{D}\cdot \Delta f_\textrm{max} \quad$ then enforce:
\vspace{-3mm}
\begin{multline} 
\textrm{T}_{1,2} \int_{0}^{\textrm{t}_2}\int_{0}^{t}\left[\Delta f_1(\tau)-\Delta f_2(\tau)\right] \textrm{d}\tau \ \textrm{d}t
\hspace*{2mm} \leq \hspace*{2mm} \\[-3pt]
2H_1\Delta f_\textrm{max} +
\frac{R_1}{\textrm{T}_\textrm{g}}\frac{(\textrm{t}_2)^2}{2} + 
\textrm{D}_1\textrm{P}^\textrm{D}_1 \int_{0}^{\textrm{t}_2}\Delta f_1(t) \ \textrm{d}t
\end{multline}
\vspace{-6mm}
\begin{multline} 
P^\textrm{L}_2\cdot t_2 
\hspace*{2mm} \leq \hspace*{2mm}
2H_2\Delta f_\textrm{max} +
\frac{R_2}{\textrm{T}_\textrm{g}}\frac{(\textrm{t}_2)^2}{2} + 
\textrm{D}_2\textrm{P}^\textrm{D}_2 \int_{0}^{\textrm{t}_2}\Delta f_2(t) \ \textrm{d}t \\[-2pt]
+ \textrm{T}_{1,2} \int_{0}^{\textrm{t}_2}\int_{0}^{t}\left[\Delta f_1(\tau)-\Delta f_2(\tau)\right] \textrm{d}\tau \ \textrm{d}t
\end{multline}
\begin{equation*}
\cdots
\end{equation*}
Four time intervals have been used in this paper for the discretisation of the nadir constraints, as it was verified that this number provided an acceptable performance. The time intervals used here correspond to the ranges $t \in [2.5\textrm{s}, 5\textrm{s}, 7.5\textrm{s}, 10\textrm{s}]$, as the nadir must happen by time `$\textrm{T}_g$' (which in Great Britain is of 10s). Using too few segments would imply that the energy lost from the outage is overestimated, as can be seen in the constraints above.

While only PFR has been considered for the sake of simplicity (i.e.~a frequency response service that must be delivered by `$\textrm{T}_g$' seconds after a fault), note that any other frequency response service with different speed and even an activation delay can be directly included in these energy-based nadir constraints. As an example, a fast service such as Enhanced Frequency Response in Great Britain (which must be delivered by 1s after a fault, typically provided by battery storage \cite{LuisEFR}) with a certain activation delay could be included in the above constraint by simply adding the following term: 
\begin{equation}
    \frac{\textrm{EFR}}{1\textrm{s}}\cdot\frac{(\textrm{t}^* - \textrm{t}_\textrm{delay})^2}{2}
\end{equation}

The proposed constraints allow not only to consider the contribution of load damping to support the nadir, neglected in many works due to the mathematical complexity it introduces, but notably also allow to consider the power sharing through transmission corridors between different regions of the power system.

\subsection{Quasi-steady-state constraint}

Regarding the quasi-steady-state frequency in multi-region systems, this magnitude does not in practice show distinct values across the network, since the inter-area oscillations are attenuated by devices such as Power System Stabilizers \cite{KundurBook} or even appropriately controlled wind farms \cite{InterAreaOscUS}.
While the action of these oscillation-damping devices is limited before the frequency nadir (which would happen in a sub-10seconds scale in GB) and therefore can be neglected, the areas are shown to swing back together by the q-s-s in the reported studies. Therefore, the q-s-s constraint aggregates the PFR and load damping in all regions to bring frequency back to $\Delta f_\textrm{max}^\textrm{ss}$:
\begin{equation}
    \sum_i^2 R_i \geq P^\textrm{L} - \Delta f_\textrm{max}^\textrm{ss} \sum_i^2 \textrm{D}_i\cdot \textrm{P}^\textrm{D}_i
\end{equation}

\section{Conditions for frequency security: \textit{N} regions} \label{ConditionsNregions}

It has been demonstrated in Section \ref{Freq_i_N_area} that the post-fault frequency evolution in any region of an $N$-region system is that of the COI plus certain inter-area oscillations, as described by eq.~(\ref{SolutionFreq_N_areas}). Since (\ref{SolutionFreq_N_areas}) is a generalisation of the two-region case described by eq.~(\ref{SolutionFreq_2areas}), which also behaves as the COI frequency plus some oscillations, the same principles and procedures discussed in Section~\ref{Conditions2Region} to obtain RoCoF and nadir constraints for a two-region system can be directly applied to the general $N$-region case.

\section{Conclusion} \label{sectionConclusion}

This paper proposes closed-form constraints for respecting frequency limits (i.e.~maximum RoCoF and nadir) in a power system that exhibits distinct frequencies in different parts of the network. Using a mathematical model for post-fault frequency dynamics in multi-region systems, it has been demonstrated that the post-fault frequency evolution in any given region is equal to the frequency of the Centre Of Inertia plus certain inter-area oscillations. 
Using this result, conditions for regional frequency security have been deduced for the first time, through a mixed analytical-numerical approach: while it is not possible to obtain purely analytical constraints, the proposed methodology uses analytical techniques to go as far as they allow, while the procedure is completed with numerical techniques on simulation samples.

Part II of this paper demonstrates the applicability of the proposed constraints in scheduling problems. 
Several case studies are run using a Stochastic Unit Commitment model, to understand the implications of inter-area frequency oscillations in the need for procuring ancillary services in each region.

\ifCLASSOPTIONcaptionsoff
  \newpage
\fi

\bibliographystyle{IEEEtran} 
\bibliography{Luis_PhD}

\begin{thebibliography}{10}
\providecommand{\url}[1]{#1}
\csname url@samestyle\endcsname
\providecommand{\newblock}{\relax}
\providecommand{\bibinfo}[2]{#2}
\providecommand{\BIBentrySTDinterwordspacing}{\spaceskip=0pt\relax}
\providecommand{\BIBentryALTinterwordstretchfactor}{4}
\providecommand{\BIBentryALTinterwordspacing}{\spaceskip=\fontdimen2\font plus
\BIBentryALTinterwordstretchfactor\fontdimen3\font minus
  \fontdimen4\font\relax}
\providecommand{\BIBforeignlanguage}[2]{{%
\expandafter\ifx\csname l@#1\endcsname\relax
\typeout{** WARNING: IEEEtran.bst: No hyphenation pattern has been}%
\typeout{** loaded for the language `#1'. Using the pattern for}%
\typeout{** the default language instead.}%
\else
\language=\csname l@#1\endcsname
\fi
#2}}
\providecommand{\BIBdecl}{\relax}
\BIBdecl

\bibitem{KundurBook}
P.~Kundur, \emph{Power System Stability and Control}, 1st~ed.\hskip 1em plus
  0.5em minus 0.4em\relax McGraw-Hill Education, 1994.

\bibitem{InmaMultiArea}
I.~Mart\'{i}nez-Sanz, B.~Chaudhuri, A.~Junyent-Ferr\'{e}, V.~Trovato, and
  G.~Strbac, ``Distributed vs. concentrated rapid frequency response provision
  in future {G}reat {B}ritain system,'' in \emph{2016 IEEE Power and Energy
  Society General Meeting}, Boston (USA).

\bibitem{DouglasWilsonSmartFreq}
P.~Wall, N.~Shams, V.~Terzija, V.~Hamidi, C.~Grant, D.~Wilson, S.~Norris,
  K.~Maleka, C.~Booth, Q.~Hong, and A.~Roscoe, ``Smart frequency control for
  the future {GB} power system,'' in \emph{2016 IEEE PES Innovative Smart Grid
  Technologies Conference Europe (ISGT-Europe)}, Ljubljana (Slovenia).

\bibitem{OverbyeLocationalInertia}
T.~Xu, W.~Jang, and T.~J. Overbye, ``Investigation of inertia's locational
  impacts on primary frequency response using large-scale synthetic network
  models,'' in \emph{2017 IEEE Power and Energy Conference at Illinois},
  Champaign (USA).

\bibitem{GoranAndersson1}
A.~Ulbig, T.~S. Borsche, and G.~Andersson, ``Impact of low rotational inertia
  on power system stability and operation,'' \emph{IFAC Proceedings Volumes},
  vol.~47, no.~3, pp. 7290--7297, 2014.

\bibitem{NatureDatabaseFreq}
L.~{Rydin Gorjao}, R.~Jumar, H.~Maass, V.~Hagenmeyer, G.~{Cigdem Yalcin},
  J.~Kruse, M.~Timme, C.~Beck, D.~Witthaut, and B.~Schafer, ``Open database
  analysis of scaling and spatio-temporal properties of power grid
  frequencies,'' \emph{Nature Communications}, vol.~11, p. 6362, 2020.

\bibitem{TimeDelaysMultiFreq}
P.~C. Bottcher, A.~Otto, S.~Kettemann, and C.~Agert, ``Time delay effects in
  the control of synchronous electricity grids,'' \emph{Chaos}, vol.~30, p.
  013122, 2020.

\bibitem{AtiaInertiaHeterogeneity}
A.~Adrees, J.~V. Milanovi\'{c}, and P.~Mancarella, ``Effect of inertia
  heterogeneity on frequency dynamics of low-inertia power systems,'' \emph{IET
  Generation, Transmission \& Distribution}, vol.~13, no.~14, pp. 2951--2958,
  2019.

\bibitem{DorflerOptimalPlacement}
B.~K. Poolla, S.~Bolognani, and F.~D\"{o}rfler, ``Optimal placement of virtual
  inertia in power grids,'' \emph{IEEE Transactions on Automatic Control},
  2017.

\bibitem{CostasVournasMultiArea}
C.~D. {Vournas} and J.~C. {Mantzaris}, ``Quasi-steady-state modeling of
  interarea oscillations,'' in \emph{2009 IEEE Bucharest PowerTech}, 2009.

\bibitem{FrequencyDivider}
F.~Milano and A.~Ortega, ``Frequency divider,'' \emph{IEEE Transactions on
  Power Systems}, vol.~32, no.~2, pp. 1493--1501, 2017.

\bibitem{RasoulPLossEstimate}
R.~Azizipanah-Abarghooee, M.~Malekpour, M.~Paolone, and V.~Terzija, ``A new
  approach to the online estimation of the loss of generation size in power
  systems,'' \emph{IEEE Transactions on Power Systems}, vol.~34, no.~3, pp.
  2103--2113, 2019.

\bibitem{SpatioTemporalMultiArea}
N.~Ma and D.~Wang, ``Extracting spatial-temporal characteristics of frequency
  dynamic in large-scale power grids,'' \emph{IEEE Transactions on Power
  Systems}, vol.~34, no.~4, pp. 2654--2662, 2019.

\bibitem{ExpositoBook}
A.~Gomez-Exp\'{o}sito, A.~J. Conejo, and C.~Ca\~{n}izares, \emph{Electric
  Energy Systems: Analysis and Operation}, 2nd~ed.\hskip 1em plus 0.5em minus
  0.4em\relax CRC Press, 2018.

\bibitem{CoherencyAnalysis}
J.~Chow, \emph{Power System Coherency and Model Reduction}, 1st~ed.\hskip 1em
  plus 0.5em minus 0.4em\relax Springer-Verlag New York, 2013.

\bibitem{FeiStochastic}
F.~Teng, V.~Trovato, and G.~Strbac, ``Stochastic scheduling with
  inertia-dependent fast frequency response requirements,'' \emph{IEEE
  Transactions on Power Systems}, vol.~31, no.~2, pp. 1557--1566, 2016.

\bibitem{LuisEFR}
L.~Badesa, F.~Teng, and G.~Strbac, ``Simultaneous scheduling of multiple
  frequency services in stochastic unit commitment,'' \emph{IEEE Transactions
  on Power Systems}, vol.~34, no.~5, pp. 3858--3868, 2019.

\bibitem{InterAreaOscUS}
Y.~Liu, J.~R. Gracia, T.~J. King, and Y.~Liu, ``Frequency regulation and
  oscillation damping contributions of variable-speed wind generators in the
  {U.S.} eastern interconnection ({EI}),'' \emph{IEEE Transactions on
  Sustainable Energy}, 2015.

\end{thebibliography}

\vskip -2\baselineskip plus -1fil

\begin{IEEEbiographynophoto}{Luis Badesa}
(S'14-M'20) received the Ph.D. degree in Electrical Engineering from Imperial College London, U.K., in 2020. He is currently a Research Associate within the Control \& Power research group at Imperial College London. His research interests lie in modelling and optimisation for low-carbon power grids.
\end{IEEEbiographynophoto}

\vskip -2\baselineskip plus -1fil

\begin{IEEEbiographynophoto}{Fei Teng}
(M'15) received the Ph.D. degree in Electrical Engineering from Imperial College London, U.K., in 2015. Currently, he is a Lecturer in the Department of Electrical and Electronic Engineering, Imperial College London, U.K. His research focuses on scheduling and market design for low-inertia power system, cyber-resilient energy system operation and control, and objective-based data analytics for future energy systems.
\end{IEEEbiographynophoto}

\vskip -2\baselineskip plus -1fil

\begin{IEEEbiographynophoto}{Goran Strbac}
(M'95) is Professor of Electrical Energy Systems at Imperial College London, U.K. His current research is focused on the optimization of operation and investment of low-carbon energy systems, energy infrastructure reliability and future energy markets. 
\end{IEEEbiographynophoto}

\end{document}